\documentclass[final]{amsart}

\usepackage{amsmath}
\usepackage{amssymb}
\usepackage{amsfonts}
\usepackage{amsthm}
\usepackage{enumerate}
\usepackage{xcolor}

\begin{document}

\newtheorem{theorem}{Theorem}[section]
\newtheorem{lemma}[theorem]{Lemma}
\newtheorem{proposition}[theorem]{Proposition}
\newtheorem{corollary}[theorem]{Corollary}

\theoremstyle{definition}
\newtheorem{definition}[theorem]{Definition}
\newtheorem{example}[theorem]{Example}
\newtheorem{examples}[theorem]{Examples}
\newtheorem{question}[theorem]{Open Question}
\newtheorem{xca}[theorem]{Exercise}
\newtheorem{remark}[theorem]{Remark}
\newtheorem{remarks}[theorem]{Remarks}

\numberwithin{equation}{section}


\def\thm#1{\begin{theorem}\label{thm #1}}
\def\endthm{\end{theorem}}
\def\thmref#1{Theorem \ref{thm #1}}
\def\cor#1{\begin{corollary}\label{cor #1}}
\def\sect#1#2{\section {#2}\label{sect #1}}
\def\sectref#1{Section \ref{sect #1}}
\def\eqn#1{\begin{equation}\label{#1}}
\def\endeqn{\end{equation}}
\def\eqnref#1{(\ref{#1})}


\def\thm#1{\begin{theorem}\label{thm #1}}
\def\endthm{\end{theorem}}
\def\thmref#1{Theorem \ref{thm #1}}
\def\cor#1{\begin{corollary}\label{cor #1}}
\def\corref#1{Corollary \ref{cor #1}}
\def\endcor{\end{corollary}}
\def\prop#1{\begin{proposition}\label{prop #1}}
\def\endprop{\end{proposition}}
\def\propref#1{Proposition \ref{prop #1}}
\def\lem#1{\begin{lemma}\label{lem #1}}
\def\endlem{end{lemma}}
\def\lemref#1{Lemma \ref{lem #1}}
\def\ex#1{\begin{example}\label{ex #1}}
\def\endex{\end{example}}
\def\exs#1{\begin{examples}\label{ex #1}}
\def\endexs{\end{examples}}
\def\exref#1{Example \ref{ex #1}}
\def\exrefs#1{Examples \ref{ex #1}}
\def\rem#1{\begin{remark}\label{rem #1}}
\def\endrem{\end{remark}}
\def\rems#1{\begin{remarks}\label{rem #1}}
\def\endrems{\end{remarks}}
\def\remref#1{Remark \ref{rem #1}}
\def\defn#1{\begin{definition}\label{defn #1}}
\def\enddefn{\end{definition}}
\def\defnref#1{Definition \ref{defn #1}}
\def\quest#1{\begin{question} \label{quest #1}}
\def\endquest{\end{question}}
\def\questref#1{Open Question \ref{quest #1}}
\def\eqn#1{\begin{equation}\label{#1}}
\def\endeqn{\end{equation}}
\def\eqnref#1{(\ref{#1})}
\def\sect#1#2{\section {#2}\label{sect #1}}
\def\sectref#1{Section \ref{sect #1}}
\def\mp#1{\marginpar{\footnotesize #1}}


\newcommand{\av}[1]{\left|#1\right|}
\newcommand{\nm}[1]{\left\|#1\right\|}
\newcommand{\cb}[1]{\left\{#1\right\}}
\newcommand{\bk}[1]{\left(#1\right)}


\def\R{{\mathbb R}}
\def\Z{{\mathbb Z}}
\def\C{{\mathbb C}}
\def\N{{\mathbb N}}
\def\T{{\mathbb T}}
\def\l{\lambda}
\def\o{\omega}
\def\re{\operatorname{Re}}
\def\sub{\subseteq}
\def\la{\langle}
\def\ra{\rangle}
\def\sap{\sigma_{\mathrm{ap}}}
\def\sp{\sigma_{\mathrm{p}}}
\def\B{\mathcal{B}}
\def\F{\mathcal{F}}
\def\D{\mathcal D}
\def\x{\mathbf{x}}
\def\a{\alpha}
\def\ep{\varepsilon}
\def\wh{\widehat}
\def\Re{\operatorname{Re}}
\def\rest{\!\!\upharpoonright}
\def\ran{\operatorname{Ran}}
\def\ker{\operatorname{Ker}}
\def\s{\sigma}

\large
\title{Lower bounds for unbounded operators and semigroups}

\author{Charles J. K. Batty}
\address{St.\ John's College, University of Oxford, Oxford OX1 3JP, England}
\email{charles.batty@sjc.ox.ac.uk}
\author{Felix Geyer}
\address{St.\ John's College, University of Oxford, Oxford OX1 3JP, England}
\email{fel.geyer@gmail.com}


\def\today{\number\day \space\ifcase\month\or
 January\or February\or March\or April\or May\or June\or
 July\or August\or September\or October\or November\or December\fi
 \space \number\year}
\date{\today}

\thanks{}


\begin{abstract} 
Let $A$ be an unbounded operator on a Banach space $X$.  It is sometimes useful to improve the operator $A$ by extending it to an operator $B$ on a larger Banach space $Y$ with smaller spectrum.   It would be preferable to do this with some estimates for the resolvent of $B$, and also to extend bounded operators related to $A$, for example a semigroup generated by $A$.  When $X$ is a Hilbert space, one may also want $Y$ to be Hilbert space.   Results of this type for bounded operators have been given by Arens, Read, M\"uller and Badea, and we give some extensions of their results to unbounded operators and we raise some open questions.  A related problem is to improve properties of a $C_0$-semigroup satisfying lower bounds by extending it to a $C_0$-group on a larger space or by finding left-inverses.   Results of this type for Hilbert spaces have been obtained by Louis and Wexler, and by Zwart, and we give some additional results.   \end{abstract}

\maketitle

\sect1{Introduction}

Let $U$ be a bounded operator on a Banach space $X$, and suppose that there is a constant $c>0$ such that
\eqn{1.0.0}
\|Ux\| \ge c\|x\|  \qquad (x \in X).
\endeqn
Then $U$ can be extended to a bounded invertible operator on a Banach space $Y$ which contains $X$ as a closed subspace; see \propref{6.0.0}, for example.   One may wish to preserve various properties of $U$.      For this to be most effective it is desirable that bounded operators which commute with $U$ can also be extended to $Y$.   This property is closely related to a result of Arens \cite{Are} for commutative Banach algebras.  The following theorem is a formulation of Arens's result for a bounded operator $U$, and $\{U\}'$ denotes the commutant of $\{U\}$ in $\B(X)$.    It can be proved in a similar way to Arens's result (see \thmref{3.3.1}).  A slightly weaker version with $\{U\}'$ replaced by any commutative subalgebra can be deduced from the result for Banach algebras, as in \cite{Mu88}.  

\thm{3.1.1}
Let $U \in \mathcal{B}(X)$ and assume that {\rm\eqnref{1.0.0}} holds for some $c>0$.  Then there exist a Banach space $Y \supseteq X$ and a unital isometric algebra homomorphism $\varphi : \{U\}' \to \mathcal{B}(Y)$ such that 
\begin{enumerate}[\rm(a)]
\item $\varphi(U)$ is invertible and $\|\varphi(U)^{-1}\| \le c^{-1}$,
\item For each $V \in \{U\}'$, $\varphi(V)$ is an extension of $V$.
\end{enumerate}
\endthm

In general, one cannot simultaneously extend two commuting operators each of which satisfy a lower bound so that both operators have inverses with the optimal norm (see \cite[Theorem 2.1]{Bol} and \cite[Example 2.3]{BY}).  Nevertheless Read extended Arens's result firstly in the context of commutative Banach algebras \cite{Rea84} and then for operators \cite{Rea88} in the form of \thmref{3.6.1} below.  Again a slightly weaker version for operators can be quickly deduced from the version for Banach algebras (see \cite[Chapter 2, Theorem 22]{Mu03}).   The full statement about the homomorphism $\varphi$ is not included in these references, but it can be seen from the proof in \cite{Rea88}.  

\thm{3.6.1}  Let $U \in \mathcal{B}(X)$.  There is a Banach space $Y \supseteq X$ and a unital isometric algebra homomorphism $\varphi : \{U\}' \to \mathcal{B}(Y)$, such that $\varphi(W)$ is an extension of $W$ for all $W \in \{U\}'$ and $\sigma(\varphi(U)) = \sap(U)$.
\endthm

Arens's theorem provides the optimal estimate for the norm of $\varphi(U)^{-1}$ whereas Read's theorem and its proof do not provide estimates for the resolvent of $\varphi(U)$.   The following result of Badea and M\"uller \cite[Theorem 3.1]{BM} concerns the norms of the powers of the inverse of the extension of a bounded operator.

\thm{3.2.3}
Let $c: \N \to (0,\infty)$ be a submultiplicative sequence.  A bounded operator $U \in \mathcal{B}(X)$ has an invertible extension $V$ on a Banach space $Y \supseteq X$ with $\|V^{-j}\| \le c_j$ for all $j\ge1$ if and only if
$$
\|x\| \le c_n\|x_0\| + c_{n-1}\|x_1\| + \dots + c_1\|x_{n-1}\|
$$
whenever $n \in \N$, $x,x_j \in X$ and $U^n x = x_0 + Ux_1 + \dots + U^{n-1} x_{n-1}$.   

When these conditions are satisfied, one may choose $Y$ and $V$ in such a way that the following hold:
\begin{enumerate}[\rm(a)]
\item $\|V^j\| = \|U^j\|$ and $\|V^{-j}\| \le c_j$ for all $j\ge1$, and
\item There is a unital isometric algebra homomorphism $\varphi : \{U\}' \to \mathcal{B}(Y)$ with $\varphi(U) = V$ and $\varphi(W)$ is an extension of $W$ for all $W \in \{U\}'$.
\end{enumerate}
\endthm 

In this paper we consider questions of this type in two further contexts, firstly replacing the single bounded operator $U$ by an unbounded operator $A$, and secondly replacing $U$ by a $C_0$-semigroup.   The common theme is that the operators should satisfy lower bounds, and they should be extended to operators on a larger Banach space with corresponding inverses, and preserving other properties as far as possible.  There is an elementary construction in \cite{MMN} which extends many bounded operators which are bounded below and it can be extended to unbounded operators (see \propref{6.0.0}).  However we seek to extend operators in the commutant of $A$ and/or to remove large parts of the spectrum, and to obtain estimates for the norms of associated bounded operators on the extended space, similarly to the theorems stated above.   The homomorphisms on $\{A\}'$ will be used to show that the constructions of Arens and Read can be carried through for generators of $C_0$-semigroups (Sections \ref{sect 3} amd \ref{sect 4}, respectively).    

For a $C_0$-semigroup $T$ on $X$ we consider lower bounds on the semigroup, of the form
$$
\|T(t)x\| \ge c(t)\|x\| \qquad (x \in X,t\ge0).
$$
In this context we seek an extension to a $C_0$-group $S$ on a larger space $Y$, as a natural analogue of the discrete case.  Alternatively one might look for a $C_0$-semigroup $L$ of left inverses on $Y$, or even on $X$ itself, so that $L(t)T(t)x=x$ for all $x \in X$ and $t\ge0$.  There are existing results in the literature for the case when $c(t)=1$ (see \propref{3.5.7b}) and for a more general case which is a continuous analogue of \thmref{3.2.3} (see \thmref{3.2.4}).   When $c(t)$ is a constant in $(0,1)$ we show that one can obtain $S$ with the exponential growth bound of  $\|S(-t)\|$ arbitrarily small (see \propref{3.2.13}), but in general not with $\|S(-t)\|$ bounded (see \exref{3.2.12}).   We also characterise when $-A$ is dissipative in an equivalent norm on $X$ (see \thmref{XS}).

In Sections \ref{sect 6} and \ref{sect 7} we consider questions of the same type in the context of operators on Hilbert spaces.  Most of the theorems above have versions in which both $X$ and $Y$ are Hilbert spaces, although some changes of detail are needed.  In addition the Hilbert space structure allows different approaches.   For example in \propref{3.3.5} we present a construction, via polar decomposition, of an extension of any closed, densely defined operator on Hilbert space satisfying a lower bound.  In \thmref{3.5.3} we show that many dissipative operators $A$ on Hilbert space $X$ have extensions to generators $B$ of contraction semigroups on larger Hilbert spaces with $D(B) \cap X = D(A)$.   Lower bounds for an operator semigroup on Hilbert space, and the possibility of finding a left-inverse semigroup on the same space, have already been considered in the literature relating to the Weiss conjecture on admissibility of observation operators in control theory \cite{LW}, \cite{XLY}, \cite{XS}, \cite{Zwa}.  In \thmref{3.2.14} we show how our results on extensions of semigroups with lower bounds relate to some of those results.

\sect2{Preliminaries}

In this paper, $X$, $Y$ and $Z$ will denote complex Banach spaces, and $\mathcal{B}(X)$ will denote the space of bounded linear operators on $X$.  We shall write $X \sub Y$, or $Y \supseteq X$, to mean that $X$ is a closed subspace of $Y$ with the same norm.  We shall also consider embeddings $\pi: X \to Y$.   When $\pi$ is isometric we regard $X$ as being a subspace of $Y$ by identifying $X$ with $\pi(X)$.  Occasionally we will allow embeddings which are isomorphisms, or continuous injections, and we will say this explicitly whenever it arises.

Given Banach spaces $X_n\; (n \in \N)$, we may consider their $\ell^p$-direct sum for $1 \le p <\infty$, or their $c_0$-direct sum, in the usual way.  In the special case when $X_n = X$ for all $n$, we will denote this space by $\ell_p(X)$ or $c_0(X)$.

An operator  $A$ on $X$ should be taken to be unbounded unless specified otherwise.  Thus the domain of $A$ is a subspace of $X$ and $A$ is a linear mapping into $X$.  We denote the domain, kernel, range, spectrum and resolvent set of $A$ by $D(A)$, $\ker A$, $\ran A$, $\s(A)$ and $\rho(A)$ respectively, and we put $R(\l,A) = (\l-A)^{-1}$ for $\l \in \rho(A)$. Recall that a \textit{pseudo-resolvent} on $X$ is a function $R : \Omega \sub \C \to \mathcal{B}(X)$ which satisfies the resolvent identity 
$$
R(\l) - R(\mu) = (\mu-\l) R(\l)R(\mu)   \qquad (\l,\mu \in \Omega).
$$
Then $\ker R(\l)$ and $\ran R(\l)$ are both independent of $\l$.  Moreover, $R$ is the resolvent of some operator $A$ on $X$ if and only if $\ker R(\l) = \{0\}$, and then $D(A) = \ran R(\l)$.  If $R(\l) = R(\l,A)$ for some $\l \in \Omega$, then $\Omega \sub \rho(A)$ and $R_\l = R(\l,A)$ for all $\l\in\Omega$  \cite[Section VIII.4]{Yos}. 

We will say that $A$ is {\it bounded below} if there exists $c>0$ such that
\eqn{2.1}
\|Ax\| \ge c\|x\|  \qquad (x \in D(A)).
\endeqn

We denote by  $\{A\}'$ the commutator algebra of $A$ in $\mathcal{B}(X)$ defined in the following way:
$$
\{A\}' = \cb{U \in \mathcal{B}(X):  \text{for all $x \in D(A)$, $Ux \in D(A)$ and $AUx=UAx$}}.
$$
If $\l\in\rho(A)$, then $\{A\}' = \{R(\l,A)\}'$ \cite[Proposition B.7]{ABHN}.

Let $A$ be an operator on $X$ and $Y$ be a Banach space with $X \sub Y$.  An operator $B$ on $Y$ is an {\it extension} of $A$ if $D(A) \sub D(B)$ and $Bx = Ax$ for all $x \in D(A)$.  We shall say that $B$ is an {\it outer extension} of $A$ if $D(A) = D(B) \cap X$ and $Bx=Ax$ for all $x \in D(A)$.

Let $T$ be a $C_0$-semigroup on $X$, so that $T : [0,\infty) \to \mathcal{B}(X)$ is continuous in the strong operator topology, with $T(0) = I$ and $T(s)T(t) = T(s+t) \; (s,t\ge0)$, and let $A$ be the generator of $T$.  We refer the reader to \cite[Section 3.1]{ABHN} or \cite[Chapter 2]{EN00} for standard properties of $C_0$-semigroups and their generators.   It is easily seen that
$$
\{A\}' = \cb{U \in \mathcal{B}(X):  \text{$T(t)U = UT(t)$ for all $t\ge0$}}.
$$

 If $S$ is a $C_0$-semigroup on $Y \supseteq X$ and $T$ is a $C_0$-semigroup on $X$, we shall say that $S$ is an \textit{extension} of $T$ if $S(t)x = T(t)x$ for all $x \in X$ and $t\ge 0$. 

We note the following elementary facts.

\prop{2.0.0}  Let $A$ be the generator of a $C_0$-semigroup $T$ on a Banach space $X$ and $B$ be the generator of a $C_0$-semigroup $S$ on a Banach space $Y \supseteq X$.
\begin{enumerate}[\rm(a)]
\item The following are equivalent:
\begin{enumerate}[\rm(i)]
\item $B$ is an extension of $A$,
\item $S(t)x = T(t)x$ for all $x \in X, \, t\ge0$,
\item $B$ is an outer extension of $A$.
\end{enumerate}
\item The following are equivalent:
\begin{enumerate}[\rm(i)]
\item $B$ is an extension of $-A$,
\item $S(t)T(t)x = x$ for all $x \in X, \, t\ge0$.
\end{enumerate}
\end{enumerate}
\endprop

\proof (a): Suppose that $B$ is an extension of $A$, and let $x \in D(A), \, t\ge0$.  Then, for $0\le s \le t$,
$$
\frac{d}{ds} S(t-s)T(s)x = - S(t-s)BT(s)x +  S(t-s)AT(s)x = 0.
$$
Hence $S(t)x = T(t)x$ for all $x\in D(A)$ and then for all $x \in X$ by density of $D(A)$ in $X$.  Thus (i)$\Rightarrow$(ii).  The proofs of (ii)$\Rightarrow$(iii)$\Rightarrow$(i) are very simple.

\noindent
(b):  The proof that (i)$\Rightarrow$(ii) is similar to (a), showing that $S(t)T(t)x$ has derivative $0$ when $x \in D(A)$, and then using the density of $D(A)$.  Now assume (ii), and let $x \in D(A)$ and $\tau>0$.  Then for $0< t < \tau$,
$$
t^{-1} (S(t) T(\tau)x - T(\tau) x) = t^{-1} (T(\tau-t)x - T(\tau)x) \to - T(\tau)Ax
$$
as $t\to0+$. Hence $T(\tau)x \in D(B)$ and $B T(\tau)x = - T(\tau)Ax$.  Letting $\tau \to 0+$, and using that $B$ is closed, it follows that $x \in D(B)$ and $Bx = -Ax$.
\endproof  

\prop{2.2.0}
Let $A$ be the generator of a $C_0$-semigroup $T$ on a Banach space $X$, let $Y$ be a Banach space with $X \sub Y$, and let $\varphi : \{A\}' \to \mathcal{B}(Y)$ be an isometric unital algebra homomorphism such that $\varphi(U)$ extends $U$ for all $U \in \{A\}'$.  Then there exists a Banach space $Y_0$ with $X \sub Y_0 \sub Y$ and a $C_0$-semigroup $S$ on $Y_0$, with generator $B$,  such that 
\begin{enumerate}[\rm(a)]
\item  $S$ extends $T$,
\item  $Y_0$ is invariant under $\varphi(U)$ for all $U \in \{A\}'$, and hence the map $\varphi_0 : U \mapsto \varphi(U)|_{Y_0}$ is an isometric unital homomorphism of $\{A\}'$ into $\mathcal{B}(Y_0)$,
\item  $\s(\varphi_0(U)) \sub \s(\varphi(U))$ for all $U \in \{A\}'$,
\item  $\s(B) \sub \s(A)$ and $R(\l,B) = \varphi_0(R(\l,A))$ for all $\l \in \rho(A)$.
\end{enumerate}
\endprop

\proof
For $t\ge0$, $T(t) \in \{A\}'$.  Let $\tilde S(t) = \varphi(T(t))$.  Then $\tilde S$ satisfies the semigroup property and is locally bounded.  Let
$$
Y_0 = \left\{ y \in Y : \lim_{t\to0+} \big\|\tilde S(t)y - y\big\| = 0 \right\}.
$$
Then $Y_0$ is a closed $\tilde S$-invariant subspace of $Y$, containing $X$.    For $U \in \{A\}'$, $\varphi(U)$ commutes with $\tilde S(t)$, so $Y_0$ is invariant under $\varphi(U)$ and under $R(\l,\varphi(U))$ for all $\l \in\rho(\varphi(U))$.  Define $\varphi_0(U) = \varphi(U)|_{Y_0}$, and $S(t) = \varphi_0(T(t))$.  Then $\varphi_0$ is a  homomorphism and $S$ is a $C_0$-semigroup on $Y_0$.  Moreover $\l \in \rho(\varphi_0(U))$ and $R(\l, \varphi_0(U)) = R(\l,\varphi(U))|_{Y_0}$.

Now, take $\l \in \rho(A)$ and let $R_\l = \varphi_0(R(\l,A))$.  Now
\begin{align*}
\left\| t^{-1}(e^{-\l t}S(t)-I) R_\l^2 + R_\l \right\|
&= \left\| t^{-1}(e^{-\l t}T(t)-I) R(\l,A)^2 + R(\l,A) \right\|  \\
&= \left\| \frac{1}{t} \int_0^t \left(I - e^{-\l s} T(s)\right) R(\l,A) \, ds \right\| \\
&\to 0 \qquad \text{as $t\to 0+$}.
\end{align*}
This shows that $R_\l = (\l - B) R_\l^2$ on $Y_0$, and hence $y = (\l-B) R_\l y$ for $y$ in the range of $R_\l$ and then for $y \in Y_0$, since $B$ is closed.  Hence $R(\l,B) = \varphi_0(R(\l,A))$ for all $\l \in \rho(A)$.   
\endproof

In the notation of the proof above, for sufficiently large  real $\l$, 
\begin{align*}
\|\tilde S(t)\varphi(R(\l,A)) - \varphi(R(\l,A))\| &= \|T(t)R(\l,A) - R(\l,A)\| \\
&\le \int_0^t \|e^{-\l s}T(s)\| \, ds \to 0
\end{align*}
as $t\to0+$.  It follows that $Y_0$ contains the range of $\varphi(R(\l,A))$.  On the other hand, $Y_0$ is the closure of the range of $R(\l,B)$ which is contained in the range of $\varphi(R(\l,A))$.   So $Y_0$ is the closure of the range of $\varphi(R(\l,A))$ for any $\l\in\rho(A)$.

\sect3{Lower bounds for unbounded operators}

There is an elementary construction which gives an invertible extension of a closed operator $A$ which is bounded below and has complemented range, but it lacks the homomorphism of the commutant. For bounded operators this construction has been given in \cite[Theorem 3]{MMN} (see also \cite{BP95} and \cite{Di00}). Our context is different and we shall need some additional properties of the construction which are not explicit in those references, so we give the details in \propref{6.0.0}. 
  
When $A$ is closed and bounded below, $\ran A$ is a closed subspace of $X$.  We assume also that $\ran A$ is complemented in $X$, so that there is a closed subspace $Y$ of $X$ such that $X = \ran A \oplus Y$.  Equivalently there is a left-inverse operator $L \in \B(X)$ such that $L$ maps $X$ into $D(A)$ and $LAx = x$ for all $x \in X$. Then let $Z$ be the $\ell_p$-direct sum, or the $c_0$-direct sum, of $X$ and countably many copies of $Y$.  Thus $Z$ consists of appropriate sequences
$$
z = (x, y_1, y_2, \dots )
$$
where  $x \in X$ and $y_n \in Y$, and
$$
\|z\| =  \big\| (\|x\|, \|y_1\|, \|y_2\|, \dots) \big\|_{\ell_p}.
$$
Then $X \sub Z$, via the isometric embedding $\pi: x \mapsto (x,0,0,\dots)$.

\prop{6.0.0}  Let $A$ be a closed operator on a Banach space $X$, satisfying {\rm\eqnref{2.1}}, and assume that $\ran A$ is complemented in $X$.  Let $Y$ and $Z$ be as above, and define an operator $B$ on $Z$ by
\begin{gather*}
D(B) = \{(x,y_1,y_2,\dots): x \in D(A), \, y_j \in Y\},  \\
 B(x, y_1, y_2, \dots) =  (Ax + c y_1, c y_2, c y_3, \dots).
\end{gather*}
Then $B$ has the following properties:
\begin{enumerate}[\rm(a)]
\item \label{3.1a} $B$ is an outer extension of $A$.
\item \label{3.1b} $B$ is invertible.
\item \label{3.1c} $\s(B) \sub \s(A)$.
\item \label{3.1d}If $A$ generates a $C_0$-semigroup on $X$, then $B$ generates a $C_0$-semigroup on $Z$.
\item \label{3.1e}If $A$ is bounded, then $B$ is bounded.
\item \label{3.1f}If $|\l|=c$, the following hold:
\begin{enumerate}[\rm(i)]
\item $\ker (\l-A) = \ker (\l-B)$, 
\item $\ran (\l-A) = \ran (\l-B) \cap X$,
\item $\ran(\l-B)$ is dense in $Z$ if $\ran(\l-A)$ is dense in $X$.
\end{enumerate}
\end{enumerate}
Moreover $B$ is a minimal invertible extension of $A$ in the sense that there is no proper closed subspace of $Z$ which contains $X$ and is invariant under $B^{-1}$.
\endprop

\proof
By replacing $A$ by $c^{-1}A$, we may assume that $c=1$. If $A$ is invertible then $Y = \{0\}$ and all properties are trivial.

Properties (\ref{3.1a}) and (\ref{3.1e}) are immediate.  For (\ref{3.1b}), $B^{-1}$ is given by
$$
B^{-1}(Ax+y_0, y_1, y_2,  \dots) =  (x, y_0, y_1, y_2, \dots).
$$

The assumption \eqnref{2.1} implies that any point $\l$ with $|\l| < 1$ is not in $\sap(A)$ and hence is not in the boundary of $\s(A)$.  If $A$ is not invertible, it follows that any point $\l \in \rho(A)$ satisfies $|\l| > 1$.  Then $R(\l,B)$ is given by
\begin{multline*}
R(\l,B)(x,y_1,y_2,\dots) \\
= \left( R(\l,A) \Big(x + \sum_{n=0}^\infty \l^{-(n+1)}y_n \Big),  \sum_{n=0}^\infty \l^{-(n+1)}y_{n+1},  \sum_{n=0}^\infty \l^{-(n+1)}y_{n+2}, \dots \right).
\end{multline*}
This establishes (\ref{3.1c}).

For (\ref{3.1f}), consider $\l$ with $|\l| = 1$. If $z = (x,y_1,y_2, \dots) \in D(B)$,  and $(\l-B)z \in \pi(X)$ then $|y_n| = |y_{n+1}|$ for all $n\ge1$.  This implies that $y_n = 0$, $z = \pi(x)$ and $(\l-B)z = \pi( (\l-A)x)$.  In particular, if $(\l-B)z=0$ then $(\l-A)x=0$.  For any vector in $Z$ of the form
$$
(x,y_1,y_2,\dots,y_k,0,0,\dots),
$$
let $y'_n=0$ for $n>k$, and define $y'_n \in Y$ recursively for $n=k,k-1,\dots,1$ by $\l y'_n = y'_{n+1} + y_n$.  If $\ran(\l-A)$ is dense we may then choose $x' \in D(A)$ such that $\|\l x' - Ax' - y'_1 - x\|$ is arbitrarily small.  Hence $\ran(\l-B)$ is dense in $Z$.

If $A$ generates a $C_0$-semigroup on $X$, then the operator $z \mapsto (Ax,0,0,\dots)$ generates a $C_0$-semigroup on $Z$.  Since $B$ is a bounded perturbation of this operator, $B$ also generates a $C_0$-semigroup. 

Finally, the span of the union of $\{B^{-k}(X) : k\ge0\}$ contains all vectors $(x,y_1,y_2,\dots)$ where $y_n = 0$ for all except finitely many $n$.  These vectors are dense in $Z$, and the minimality follows.
\endproof

The following is a version of \thmref{3.1.1} for semigroup generators.  The proof here is an adaptation of Arens's proof in the context of Banach algebras.

\thm{3.3.1}
Let $A$ be the generator of a $C_0$-semigroup on $X$, and assume that $A$ satisfies {\rm\eqnref{2.1}}
for some $c>0$.  Then there exist a Banach space $Y \supseteq X$ and an operator $B$ on $Y$ with the following properties:
\begin{enumerate}[\rm(a)]
\item\label{3.2a} $B$ is the generator of a $C_0$-semigroup on $Y$,
\item\label{3.2b} $B$ is an outer extension of $A$,
\item\label{3.2c} $B$ is invertible with $\|B^{-1}\| \le c^{-1}$, and $\sigma(B) \subseteq \sigma(A)$,
\item\label{3.2d} There is a unital isometric algebra homomorphism $\varphi : \{A\}' \to \mathcal{B}(Y)$
such that $\varphi(U)$ is an extension of $U$ for all $U \in \{A\}'$ and $\varphi(R(\l,A)) = R(\l,B)$ for all $\l \in \rho(A)$.
\end{enumerate}
\endthm 

\proof
Replacing $A$ by $c^{-1}A$ we may assume that $c=1$.  Let $T$ be the $C_0$-semigroup generated by $A$.  We first lift $T$ to the space $\ell_1(X)$.

Define
$$
\tilde T(t)f = (T(t)x_n) \qquad (f=(x_n) \in \ell_1(X), t\ge0).
$$
Then $\tilde T$ is a $C_0$-semigroup on $\ell_1(X)$ and its generator $C$ is given by 
\begin{align*}
D(C) &= \big\{f=  (x_n) \in \ell_1(X):  x_n \in D(A)\; (n\in\N),  \; (Ax_n) \in \ell_1(X) \big\}, \\
Cf &= (Ax_n).
\end{align*}
Moreover $\s(C) = \s(A)$ and the resolvent of $C$ is given by
$$
R(\lambda,C)f = (R(\l,A)x_n) \qquad (\l\in \rho(A), \; f=(x_n) \in \ell_1(X)).
$$
Let $\bar R$ be the right shift on $\ell_1(X)$, and note that $\bar R \in \{C\}'$.  Let $J$ be the closure of $\{ f - C \bar Rf : f \in D(C)\}$ in $\ell_1(X)$.  Then $J$ is invariant under $\tilde T$, $\bar R$ and $R(\l,C)$.

Let $Y = \ell_1(X)/J$ and define $\pi:X \to Y$ by $\pi(x) = x\mathbf{e}_0 + J$.   Here $x\mathbf{e}_0$ is the sequence $f$ with $f(0)=x$ and $f(n)=0$ for $n\ge1$.  It is clear that $\|\pi(x)\| \le \|x\|$.  Using the triangle inequality and the lower bound \eqnref{2.1} for $A$ repeatedly, for any $f = (x_n) \in D(C)$, we have
\begin{eqnarray*}
\|x\| &\le& \|x-x_0\| + \|Ax_0\| \le \|x-x_0\| + \|x_1-Ax_0\| + \|Ax_1\| \le \cdots  \\
&\le& \|x-x_0\| + \sum_{n=1}^\infty \|x_n-Ax_{n-1}\| \\
&=& \|x\mathbf{e}_0 - (I-C\bar R)f\|_{\ell_1(X)}.
\end{eqnarray*}
Thus $\|x\| \le \|\pi(x)\|$ and $\pi$ is isometric.

Let $q : \ell_1(X) \to Y$ be the quotient map $q(f) = f + J$.  Since $J$ is invariant under $\tilde T(t)$, there are operators on $Y$ defined by
$$
S(t)q(f) = q(\tilde T(t)f) \qquad (f \in \ell_1(X), t\ge0).
$$
Then $S$ is a $C_0$-semigroup and its generator $B$ is given by $D(B) = q(D(C))$, $Bq(f) = q(Cf)$ for $f \in D(C)$.  Moreover, $\s(B) \sub \s(C) = \s(A)$ and $R(\l,B)q(f) = q(R(\l,C)f)$ for $f \in \ell_1(X)$ (see \cite[I.5.13, II.2.4 and IV.2.15]{EN00}).  For $x \in X$, we have
$$
S(t)\pi(x) = \pi(T(t)x).
$$
When we identify $X$ with its image under $\pi$, this shows that $B$ is an outer extension of $A$ (\propref{2.0.0}(a)).  

Next we show that $B$ has a bounded inverse.  Define $V$ on $Y$ by
$$
Vq(f) = q(\bar Rf) \qquad (f \in \ell_1(X)).
$$
Since $J$ is invariant under $\bar R$, the operator $V$ is well-defined and $\|V\| \le \|\bar R\|=1$.  Let $y \in D(B)$ and choose $f \in D(C)$ such that $y = q(f)$.  Since $V$ commutes with $S(t)$, $Vy \in D(B)$ and 
$$
BVy = VBy = q(C\bar Rf) = q(f) - q(f-C\bar Rf) = y.
$$
Since $B$ is closed, $D(B)$ is dense in $Y$ and $V$ is bounded, we deduce that $Vy \in D(B)$ and $BVy=y$ for all $y\in Y$, and $V$ is the inverse of $B$.

For $U \in \{A\}'$ and $f = (x_n) \in \ell_1(X)$, define
$$
\varphi(U)q(f) := q((Ux_n)).
$$
Since $J$ is invariant under the map $f \mapsto (Ux_n)$, $\varphi(U)$ is a well-defined, bounded operator on $Y$ with $\|\varphi(U)\| \le \|U\|$ and $\varphi(U)\pi(x) = \pi (Ux)$ for $x \in X$, so $\varphi(U)$ is an extension of $U$ when we identify $X$ with $\pi(X) \sub Y$.  Moreover $\varphi$ is a unital algebra homomorphism.  For $\l \in \rho(A)$ and $f \in \ell_1(X)$, we have
$$
\varphi(R(\l,A))q(f) = q((R(\l,A)x_n)) = q(R(\l,C)f) = R(\l,B)q(f).
$$
So $\varphi(R(\l,A)) = R(\l,B)$.
\endproof

It is plausible that \thmref{3.3.1} can be extended to larger classes of unbounded operators than generators of $C_0$-semigroups, but this is not straightforward.  In trying to extend the proof of \thmref{3.3.1}, one can define an operator $C$ on $\ell_1(X)$ and the corresponding space $J$.  Then the resolvent of $A$ induces operators $R_\l$ on $\ell_1(X)/J$ which form a pseudo-resolvent.  The technical problem which arises is to show that $R_\l$ is injective so that the pseudo-resolvent is the resolvent of an operator.  This can be achieved in the proof of \thmref{3.3.1}, because $C$ generates a $C_0$-semigroup which leaves $J$ invariant, so that if $f \in D(C) \cap J$ then $Cf \in J$ (\propref{2.0.0}).    A similar problem arises if one assumes that $\rho(A)$ is non-empty and applies \thmref{3.1.1} to $U := \mu^{-1} - R(\mu,A)$ where $\mu \in \rho(A)$.  

In seeking extensions it is natural to consider the class of generators of integrated semigroups, but we are able to obtain only a weak result (\propref{3.3.3}) for once integrated semigroups, with a weaker conclusion giving a continuous embedding of $X$ in $Y$.

Let $k\ge1$.  A function $T : [0,\infty) \to \B(X)$ which is continuous in the strong operator topology, is said to be a \textit{$k$-times integrated semigroup} if $T(0)=0$ and there exist $\o\ge0$ and $M\ge1$ such that 
$$
\nm{ \int_0^t T(s)x \, ds} \le Me^{\o t}\|x\| \qquad (x \in X, t\ge 0),
$$
and the (improper) integral
$$
R(\l)x := \l^k \int_0^\infty e^{-\l t} T(t)x \, dt  \qquad (x \in X, \l>\o)
$$
defines a pseudo-resolvent.  These integrals may be improper, but the integrals over $[0,\tau]$ converge in operator-norm (i.e., uniformly for $\|x\|\le1$), as $\tau\to\infty$ \cite[Remark 1.4.6]{ABHN}.  Then $T$ is {\it non-degenerate} if any one or all of the following equivalent conditions is satisfied:
\begin{enumerate}[\rm(i)]
\item $R(\l)$ is injective.
\item There is an operator $A$ such that $R(\l,A) = R(\l)$.  
\item If $x \in X$ and $T(t)x= 0$ for all $t\ge0$ then $x=0$.
\end{enumerate} 
See \cite[Proposition 3.2.9]{ABHN}.  Then $A$ is the \textit{generator} of the integrated semigroup.

When $A$ generates a non-degenerate (once) integrated semigroup $T$ on $X$ and $Y$ is a closed $T$-invariant subspace of $X$, there may exist $y \in Y \cap D(A)$ with $Ay \notin Y$, and the induced integrated semigroup on $X/Y$ is then degenerate.  This causes the proof of \thmref{3.3.1} to break down for integrated semigroups.  In the special case of once integrated semigroups we can partially avoid this obstruction as shown in \propref{3.3.3}.

Let $S$ be a (degenerate) $k$-times integrated semigroup on $X$, and $N$ be its {\it degeneration space}
$$
N := \cb{y \in X:  \text{$S(t)y =0$ for all $t\ge0$}}.
$$
This is a closed subspace of $X$, invariant under $S(t)$.  Let $\tilde S(t)$ be the operator on $X/N$ induced by $S(t)$.   It is readily verified that $\tilde S$ is a $k$-times integrated semigroup.

\lem{2.5.10}  If $S$ is a once integrated semigroup on $X$, then $\tilde S$ is a non-degenerate once integrated semigroup on $X/N$.
\end{lemma}

\proof
We have to show that $\tilde S$ is non-degenerate, that is, if $x \in X$ and $S(t)x \in N$ for all $t\ge0$, then $x \in N$.   Since $S(t)x \in N$, $S(s)S(t)x = 0$ for all $s,t\ge0$.  By \cite[Proposition 3.2.4]{ABHN},
$$
\int_s^{s+t} S(r)x \, dr = \int_0^t S(r)x \, dr.
$$
Differentiation at $t=0$ gives $S(s)x=S(0)x=0$ for all $s$, so $x\in N$. 
\endproof

\lemref{2.5.10} does not extend to $k$-times integrated semigroups, where $k\ge2$.

\ex{2.5.6}  
Let $k \ge 2$, $X = \C^2$ and
$$
S(t) = \begin{pmatrix}   0 & t^{k-1} \\ 0&0 \end{pmatrix} \qquad (t\ge0).
$$
Then $S(0)=0$ and 
$$
R(\l) = \l^k \int_0^\infty e^{-\l t} S(t) \, dt = \begin{pmatrix} 0 & (k-1)! \\ 0&0 \end{pmatrix} \qquad (\l>0).
$$
This is a pseudo-resolvent, so $S$ is a $k$-times integrated semigroup.  The degeneration space $N$ is spanned by $(1,0)$.  The induced integrated semigroup on $X/N$ is $0$, so it is degenerate.
\endex

\prop{3.3.3}  Let $A$ be the generator of a non-degenerate once integrated semigroup on $X$, and assume that {\rm\eqnref{2.1}} holds for some $c>0$.   Then there exist a Banach space $Z$, a continuous embedding $\pi : X \to Z$ and an operator $B$ on $Z$ with the following properties:
\begin{enumerate}[\rm(a)]
\item $B$ is the generator of a non-degenerate once integrated semigroup on $Z$,
\item $B$ is invertible with $\|B^{-1}\| \le c^{-1}$, and $\s(B) \sub \s(A)$,
\item If $x \in D(A)$ then $\pi(x) \in D(B)$ and $B\pi(x) = \pi(Ax)$.
\end{enumerate}
\endprop

\proof
The proof initially follows the lines of \thmref{3.3.1}.  Let $T$ be the once integrated semigroup generated by $A$.   Define $Y$, $C$, $\bar R$, $J$, $\tilde T(t)$ and $S(t)$ as in the proof of \thmref{3.3.1}.  They have the same properties as in \thmref{3.3.1}, except that $\tilde T$ and $S$ are now once integrated semigroups, $\tilde T$ is non-degenerate with generator $C$, and $S$ may be degenerate.  Nevertheless there exist bounded operators $R_\lambda \; (\l \in \rho(A))$ on $Y$ such that  $R_\l (f + J) = (R(\l,A)x_n) + J$ for $f = (x_n) \in \ell_1(X)$.  They form a pseudo-resolvent, but they are not necessarily injective.   

Let $N$ be the degeneration space of $S$, and let $Z = Y/N$.  Since $R_\l$ commutes with $S(t)$, $N$ is invariant under $R_\l$, and therefore $R_\l$ induces a bounded operator $\tilde R_\l$ on $Z$, and $\{\tilde R_\l : \l \in \rho(A)\}$ is a pseudo-resolvent.   By \lemref{2.5.10}, $S$ induces a non-degenerate once integrated semigroup $\tilde S$ on $Z$. Let $q(f) = (f+J)+N$ for $f \in \ell_1(X)$, and $\pi(x) = q(x\mathbf{e}_0)$ for $x \in X$.  Let $x\in \ker\pi$, so $S(t)(x\mathbf{e}_0 + J)=0$ for all $t\ge0$.  This means that $T(t)x \mathbf{e}_0 \in J$, and then $T(t)x=0$ for all $t\ge0$ since $y \to y\mathbf{e}_0+J$ is injective on $X$.  Since $T$ is non-degenerate, this implies that $x=0$.  Thus $\pi$ is a continuous embedding of $X$ into $Z$.  

Define $V : Z \to Z$ by
$$
Vq(f) = q(\bar Rf) \qquad (f \in \ell_1(X)).
$$ 
This is well-defined and contractive, since $f+J \in N$ implies that $\tilde T(t)f \in J$ for all $t \ge 0$, and then $\tilde T(t)\bar Rf = \bar R\tilde T(t)f \in J$, so $\bar Rf + J \in N$.   Let $B$ be the generator of $\tilde S$, take $\l>\o$, and consider $z \in D(B)$.  Then
$$
z = R(\l,B)y = \l \int_0^\infty e^{-\l t} \tilde S(t)y \, dt
$$
for some $y \in Z$.  Take $g \in \ell_1(X)$ such that $y = q(g)$, and let $f = R(\l,C)g$.  Then $z = q(f)$ and
$$
VBz = \l Vz - Vy = \l Vz - q \bar Rg = \l Vz - q \bar R(\l-C)f = q C\bar Rf = q(f).
$$
Moreover, for any $z = q(f)$ and $\l>\o$, 
\begin{align*}
R(\l,B)(\l Vz-z) &= q R(\l,C) (\l \bar Rf-f) \\
&= \l q \bar R  R(\l,C) f - q \bar RC R(\l,C)f = q \bar Rf = Vz.
\end{align*}
Hence $Vz \in D(B)$ and $BVz = z$.  So $B$ is invertible with $B^{-1} = V$.  Moreover for $f \in \ell_1(X)$ and $\l>\o$, 
\begin{align*}
q R(\l,C)f &=  q \left( \lim_{\tau\to\infty}  \left(\l \int_0^\tau e^{-\l t} \tilde T(t) f \, dt \right)  \right) \\
&= \l \lim_{\tau\to\infty} \int_0^\tau \tilde S(t) q(f) \, dt = R(\l,B)q(f).
\end{align*}
Since $\{\tilde R_\l : \l \in \rho(A)\}$ is a pseudo-resolvent, it follows that $\s(B) \sub \s(A)$.

For $x \in D(A)$, we have $x = R(\l,A) (\l x - Ax)$, so
$$
\pi(x) = \tilde R_\l \pi(\l x - Ax) = R(\l,B) \pi (\l x - Ax).
$$
It follows that $\pi(x) \in D(B)$ and $B \pi(x) = \pi(Ax)$.
\endproof

\sect4{Reducing the spectrum}

Here we adapt Read's theorem (\thmref{3.6.1}) itself to the case of the generator of a $C_0$-semigroup.    

\thm{3.6.2}  Let $A$ be the generator of a $C_0$-semigroup on $X$.   Then there exist a Banach space $Y \supseteq X$ and an operator $B$ on $Y$ with the following properties:
\begin{enumerate}[\rm(a)]
\item $B$ is the generator of a $C_0$-semigroup on $Y$, 
\item $B$ is an outer extension of $A$,
\item $\sigma(B) = \sap(A)$.
\end{enumerate}
Furthermore there is a unital isometric algebra homomorphism 
$$
\varphi : \{A\}' \to \mathcal{B}(Y)
$$
such that $\varphi(U)$ is an extension of $U$ for all $U \in \{A\}'$ and $\varphi(R(\l,A)) = R(\l,B)$ for all $\l \in \rho(A)$.
\endthm

\proof
Let $T$ be the $C_0$-semigroup generated by $A$, and take $\mu > \o_0(T)$, the growth bound of $T$, so $\mu \in \rho(A)$. By applying \thmref{3.6.1} to $R(\mu,A)$, and then applying \propref{2.2.0},  there exist a Banach space $Y \supseteq X$ and a $C_0$-semigroup $S$ on $Y$, extending $T$, such that the generator $B$ of $S$ satisfies $\s(R(\mu,B)) \sub \sap(R(\mu,A))$.  Since $R(\mu,B)$ extends $R(\mu,A)$, we have
$$
\s(R(\mu,B)) = \sap(R(\mu,A)).
$$
By a standard spectral property \cite[Proposition B.2]{ABHN},
\begin{align*}
\sap(A) &= \cb{ \l\in\C: (\mu-\l)^{-1} \in \sap(R(\mu,A))} \\
&= \cb{ \l\in\C: (\mu-\l)^{-1} \in \sigma(R(\mu,B))} = \sigma(B). \qedhere
\end{align*}
\endproof

A class of ``uniformly expansive'' operators was introduced in \cite{EH} and \cite{Hed}, and the following slightly larger class was considered in \cite{BT10}.  An operator $U \in \mathcal{B}(X)$ is said to be \textit{quasi-hyperbolic} if the following equivalent conditions are satisfied.  Here, $\mathbb{T}$ denotes the unit circle in $\C$.
\begin{enumerate}[\rm(i)]
\item There exists $n \in \N$ such that $\max(\|U^{2n}x\|,\|x\|) \ge 2 \|U^nx\|$ for all $x \in X$.
\item $\sap(U) \cap \mathbb{T} = \emptyset$.
\item There exists $c>0$ such that
\eqn{4.2}
\|(U-\l)x\| \ge c \|x\| \qquad (x \in X, \,\l \in \mathbb{T}).
\endeqn
\item $U$ is the restriction of an operator $V \in \mathcal{B}(Y)$ on a Banach space $Y \supseteq X$ with $\mathbb{T} \sub \rho(V)$ (i.e., $V$ is hyperbolic).  
\end{enumerate}
The implication (ii)$\Rightarrow$(i) was essentially shown in \cite{Hed}, and  (ii)$\Rightarrow$(iv) is a special case of \thmref{3.6.1}.  The other implications are very elementary.

Although the proof of Read's theorem is substantially simpler in the restricted form needed here than in its full generality, it does not appear to give any estimate for $\|R(\l,V)\|$.  Nevertheless it is possible to show that there exists $K$ depending only on $\|U\|$ and $c$, such that $V$ can be constructed as in (iv) and $\|R(\l,V)\| \le K$ for all $\l \in \mathbb{T}$.  To see this, we argue by contradiction.

Let $M>1$ and $c>0$.  Assume that for each $n \in \N$, there is an operator $U_n$ on a Banach space $X_n$ such that $U_n$ satisfies \eqref{4.2} holds, $\|U_n\| \le M$ and for every extension $V$ of $U_n$ with $\mathbb{T} \sub \rho(V)$, there exists $\l \in \C$ such that $\|R(\l,V)\| > n$.  Let $X$ be the $c_0$-direct sum $X = \bigoplus_{n\ge1} X_n$ and $U = \bigoplus_{n\ge1} U_n$.  Then $\|U\|\le M$ and 
$$
\|(U-\l)(x_n)\| \ge c \|(x_n)\| \qquad ((x_n) \in X, \l \in \T).
$$
By \thmref{3.6.1} there is an extension $V$ of $U$ with $\T \sub \rho(V)$.  Since $\T$ is compact there exists $K$ such that $\|R(\l,V)\| \le K$ for all $n\ge1$.  Now $X_n$ is isometrically embedded in $X$ and $V$ is an extension of $U_n$.   For each $n\ge1$, $n < \|R(\l_n,V)\| \le K$ for some $\l_n \in \T$.  This is a contradiction.

This raises the following problem.

\quest{3.6.3}
Find $K : (0,1) \times (1,\infty) \to (1,\infty)$ such that the following holds.  If $U$ is a bounded operator satisfying \eqnref{4.2}, there is an extension $V$ of $U$ with $\|V\|=\|U\|$, $\mathbb{T}\sub \rho(V)$ and
$$
\sup_{\l\in\mathbb{T}} \|R(\l,V)\| \le K(c,\|U\|).
$$
\endquest

Weighted shift operators provide a class of quasi-hyperbolic operators which are not hyperbolic.   The following gives an explicit description of extensions of weighted shifts on $\ell_1(\mathbb{Z})$ to hyperbolic operators. 

\ex{3.6.1}
Let $X = \ell_1(\Z)$ and $T_w \in \mathcal{B}(X)$ be the weighted shift given by 
$$
T_w(\mathbf{e}_n) = w_n \mathbf{e}_{n+1}  \qquad (n\in\Z),
$$
where $\mathbf{e}_n$ are the standard basis vectors and $0 < w_n \le M$ for all $n$.   We assume that
$$
w_n \ge 1 \quad (n\ge0), \qquad w_n < 1 \quad (n<0), \qquad 0 \le r^-(T_w) < 1 < i^+(T_w).
$$
Here
\begin{align*}
i^+(T_w) &= \lim_{n\to\infty} \inf_{k>0} (w_k \dots w_{k+n-1})^{1/n},  \\
r_-(T_w) &= \lim_{n\to \infty} \sup_{k<0} (w_{k-n} \dots w_{k-1})^{1/n}.
\end{align*}
For the corresponding operator on $\ell_2(\Z)$, Ridge \cite{Rid} (see also \cite[Proposition 2.5]{BT10}) showed that $\sap(T_w)$ consists of two (possibly degenerate) annuli centred at 0, one with outer radius $r^-(T_w)$,  and the other with inner radius $i^+(T_w)$, while $\s(T_w)$ is the smallest annulus containing $\sap(T_w)$.  In particular, $T_w$ is quasi-hyperbolic, but $\mathbb{T} \sub \s(T_w)$.
 
In \cite[Section 3.6.1]{Gey} an isomorphic embedding $\pi : X \to Y := \ell_1(\mathbb{Z}) \oplus \ell_1(\mathbb{Z})$, and weights $\alpha$ and $\beta$ on $\mathbb{Z}$, are constructed such that 
\eqn{4.2.5}
\pi \circ T_w = S \circ \pi
\endeqn
where $S = T_\alpha \oplus T_\beta$.  Moreover, the weights $\alpha$ and $\beta$ are chosen so that $\s(T_\alpha)$ is inside the open unit disc and $\s(T_\beta)$ is outside the closed unit disc.  In particular, $\s(S) \cap \mathbb{T}$ is empty.   Furthermore it is shown that if $c>0$ is chosen so that $T_w$ satisfies \eqnref{4.2} then
\eqn{4.3}
\|R(\l,S)\| \le \frac{e \tilde M^{2/c+1} \log \tilde M}{c} \qquad (\l \in \mathbb{T}),
\endeqn
where $\tilde M = \max(2, \|T_w\|)$.   

The embedding $\pi$ of $X$ in $Y$ here is not isometric but $\|x\|_* := \|\pi(x)\|_Y$ is a weighted $\ell_1$- norm satisfying $\|x\|_X \le \|x\|_* \le 2 \|x\|_X$.  One may vary the definition of $\pi$ so that it becomes isometric.  Then one can vary the definition of $\alpha$ and $\beta$ in order to satisfy \eqnref{4.2.5}. 
\endex

Other constructions or calculations may provide sharper estimates than \eqnref{4.3}.   On the other hand, 
it is plausible that weighted shifts will have the worst behaviour for this type of problem.  Thus we conjecture that there exist positive constants $a,b$ such that
$$
K(c,\tau) := a \tau^{b/c}
$$
has the properties required in \questref{3.6.3}.   

There is a related question for generators of quasi-hyperbolic $C_0$-semigruoups as defined in  \cite{BT10}.  A $C_0$-semigroup $T$ is quasi-hyperbolic if and only if $T(t)$ is quasi-hyperbolic for some, or equivalently all, $t>0$.

Let $A$ be the generator of a $C_0$-semigroup $T$.  When $T$ is quasi-hyperbolic, it was shown in \cite[Proposition 3.2]{BT10} that
\eqn{3.23}
\|(A-is)x\| \ge c\|x\| \qquad (x \in D(A), s \in \R)
\endeqn
for some $c>0$.  On the other hand, the full converse is not true, as \eqnref{3.23} does not imply that $T$ is quasi-hyperbolic in general.  A possible partial converse was raised as a question in \cite[Section 4]{BT10}; it remains an open question whether \eqnref{3.23} implies that there is a continuous injection of $X$ into a Banach space $Y$ such that $\pi \circ T(t) = S(t) \circ \pi$ for some hyperbolic semigroup $S$ on $Y$.

 The following question asks about another possible partial converse property of a different type.

\quest{3.6.4}
Let $A$ be the generator of a $C_0$-semigroup on a Banach space $X$, and assume that
\eqnref{3.23} holds for some $c>0$.   Does $A$ have an extension $B$ which is the generator of a $C_0$-semigroup on a Banach space $Y \supseteq X$, with $i\R \sub \rho(B)$ and $\sup_{s\in\R} \|R(is,B)\| < \infty$?
\endquest

This question has a positive answer on Hilbert spaces.  For a $C_0$-semigroup $T$ on a Hilbert space, the following properties are equivalent:
\begin{enumerate}[\rm(i)]
\item \eqnref{3.23} holds,
\item $T$ is quasi-hyperbolic,
\item $T$ is the restriction of a hyperbolic $C_0$-semigroup on a Hilbert space to a closed invariant subspace,
\item $A$ has an extension $B$ which is the generator of a $C_0$-semigroup on a Hilbert space $Y \supseteq X$, with $i\R \sub \rho(B)$ and $\sup_{s\in\R} \|R(is,B)\| < \infty$.
\end{enumerate}
The implication (i)$\Rightarrow$(ii) was shown in \cite[Corollary 3.1]{BT10}.  Read's theorem may be applied as in \cite[Theorem 2.2]{BT10} to show that (ii)$\Rightarrow$(iii).  Elementary theory of hyperbolic semigroups as in \cite[Section V.1c]{EN00} shows that  that (iii)$\Rightarrow$(iv), and (iv)$\Rightarrow$(i) is elementary.

\sect5{Lower bounds for semigroups} 

Let $T$ be a $C_0$-semigroup on a Banach space $X$, and assume that, for some $\tau>0$, there exists $c(\tau)>0$ such that
$$
\|T(\tau)x\| \ge c(\tau)\|x\| \qquad (x \in X).
$$
Using the semigroup property and local boundedness of $\|T(t)\|$, it is easy to see that 
\eqn{4.0}
\|T(t)x\| \ge c(t) \|x\| \qquad (x \in X, t\ge 0),
\endeqn
where $c : (0,\infty) \to (0,\infty)$.  The function
$$
\gamma_T(t) := (\inf \{\|T(t)x\| : x \in X, \|x\|=1\})^{-1},
$$
is submultiplicative, i.e., $\gamma_T(s+t) \le \gamma_T(s)\gamma_T(t)$, and locally bounded, and it satisfies  $\gamma_T(t) \le M e^{\a t}$ for some $M>0$ and $\a \in \R$.  Now \eqnref{4.0} holds if and only if $c(t) \ge 1/\gamma_T(t)$.

In some parts of the literature, a $C_0$-semigroup is said to be ``left-invertible'' if \eqnref{4.0} holds with $c(t)>0$, but such terminology for semigroups might be interpreted in several different ways for semigroups.    For greater precision, we shall say that $T$ {\it satisfies lower bounds} when \eqnref{4.0} holds, with $c(t)>0$.   Now we describe various notions of left-inverse for $T$.

If each operator $T(t)$ has a bounded left-inverse on $X$, then $T$ satisfies lower bounds.  When $T$ satisfies lower bounds, each $T(t)$ is injective with closed range, so there is a bounded left-inverse operator $L_0(t) : \ran T(t) \to X$ such that $L_0(t)T(t) = I$.  In many cases, $L_0(t)$ can be extended to a bounded operator $L(t)$ on $X$, and then $L(t)T(t) = I$.  However this extension process may be quite arbitrary, and $L$ may not satisfy the semigroup property on $X$.

We shall say that a $C_0$-semigroup $L$ on $X$ is a \textit{left-inverse semigroup} on $X$ for $T$ if 
$L(t)T(t)x = x$ for all $x\in X, \,t\ge0$.  More generally, a $C_0$-semigroup $L$ on a Banach space  
 $Y \supseteq X$ is a \textit{left-inverse semigroup} on $Y$ for $T$ if $L(t)T(t)x = x$ for all $x\in X, \,t\ge0$.  \propref{2.0.0} shows that $L$ is a left-inverse semigroup for $T$ if and only if  the generator $B$ of $L$ is an extension of $-A$, where $A$ is the generator of $T$.   It is known that any $C_0$-semigroup satisfying lower bounds on a Hilbert space $X$ has a left-inverse semigroup on $X$ (see \thmref{3.2.14}).   

A particular case of a left-inverse semigroup occurs when $T$ extends to a $C_0$-group $S$ on $Y \supseteq X$  in the sense that $X$ is invariant under $S(t)$ and $T(t) = S(t)|_X$, for all $t\ge0$.  In this case, we shall say that $S$ is a \textit{$C_0$-group extension} of $T$ on $Y$.  Then $L(t):=S(-t)$ is a left-inverse semigroup for $T$ on $Y$, and the generator of $S$ is an extension of $A$.  However $X$ will not be invariant under $S(-t)$ for $t > 0$ unless $T$ is itself a $C_0$-group.  We shall see in \propref{3.2.13} that every $C_0$-semigroup satisfying lower bounds has a $C_0$-group extension on a larger Banach space. 

When constructing left-inverse $C_0$-semigroups $L$ or $C_0$-group extensions $S$ for $T$,  it is desirable also to keep control on their growth.  Ideally the growth of $\|L(t)\|$ would be comparable with the growth of $\gamma_T(t)$, but this cannot be achieved precisely in general (see \propref{3.2.13} and \exref{3.2.12}).  In the case of group extensions, it would also be desirable that $\|S(t)\|$ is comparable to $\|T(t)\|$ for large $t>0$.

We now recall some known results, beginning with the case when $c(t) = 1$ for all $t$.  In the following proposition the equivalence of (i) and (ii) is a counterpart of the Lumer-Phillips theorem, from \cite[pp.419,420]{GS}, and the equivalence of (i) and (iii) is from \cite[Proposition 2.2, Theorem 3.3]{BY}.

\prop{3.5.7b}  Let $A$ be the generator of a $C_0$-semigroup $T$ on a Banach space $X$.  The following are equivalent:
\begin{enumerate}[\rm(i)]
\item $T$ is expansive, i.e., $\|T(t)x\| \ge \|x\|$ for all $x \in X, \, t\ge0$.
\item $-A$ is dissipative. 
\item There is a $C_0$-group extension $S$ of $T$  on a Banach space $Y \supseteq X$ such that $\|S(-t)\| \le 1$ for all $t\ge0$.
\end{enumerate}
\endprop

The following theorem is a compilation of known results giving a version of \propref{3.5.7b} allowing arbitrary (submultiplicative) growth of the left-inverses.  The property (\ref{5.2a}) of \thmref{3.2.4} implies \eqnref{4.0} for $c(t) = \gamma(t)^{-1}$, by putting $n=1$ and $t_1=0$.  The converse (that \eqnref{4.0} implies (\ref{5.2a})) holds when $c(t)=1$, or $c(t) = e^{-\a t}$ (see \propref{3.5.7b}), but not in general (see \exref{3.2.12}).   The implication (\ref{5.2a})$\Rightarrow$(\ref{5.2b}) was noted in \cite[Remarks 3.2(iv)]{BM}, where details of the discrete counterpart were given.   For (\ref{5.2c})$\Rightarrow$(\ref{5.2a}), applying $L(t)$ to \eqnref{5.3} gives
$$
x = L(t-t_1) x_1 + \dots + L(t-t_n) x_n,
$$
and (\ref{5.2a}) follows.

\thm{3.2.4}
Let $\gamma : \R_+ \to (0,\infty)$ be a submultiplicative function.  Let $T$ be a $C_0$-semigroup on $X$.  The following properties are equivalent:
\begin{enumerate}[\rm(a)]
\item\label{5.2a}  Whenever $x, x_1,\dots, x_n \in X$, $0 \le t_1 \le t_2 \le \dots \le t_n \le t$, and 
\eqn{5.3}
T(t) x = T(t_1) x_1 + \dots + T(t_n) x_n,
\endeqn
we have
\eqn{5.2}
\|x\| \le  \gamma(t-t_1)\|x_1\| + \dots + \gamma(t-t_n) \|x_n\|.
\endeqn
\item\label{5.2b}  There is a $C_0$-group extension  $S$ of $T$ on a Banach space $Y \supseteq X$ such that $\|S(-t)\| \le \gamma(t)$ for all $t\ge0$.
\item\label{5.2c}  There is a left-inverse semigroup $L$ for $T$ on a Banach space $Y \supseteq X$ such that $\|L(t)\| \le \gamma(t)$ for all $t\ge0$.
\end{enumerate}
\endthm

Zwart \cite[Section 3]{Zwa} has asked whether every $C_0$-semigroup satisfying lower bounds on a Banach space $X$ has a left-inverse semigroup on $X$.  The following new result shows that any semigroup satisfying lower bounds has a $C_0$-group extension $S$ on a larger Banach space such that the exponential growth bound of $\{S(-t) : t\ge0\}$ is arbitrarily close to the exponential growth bound of $\gamma_T$.  In \exref{3.2.12} we show that it may not be possible to arrange that $\|S(-t)\| = O\left(c_T(t)^{-1}\right)$ as $t\to\infty$.

\prop{3.2.13}   Let $T$ be a $C_0$-semigroup on a Banach space $X$ satisfying lower bounds, and let $c>0$ and $\a\in\R$ be such that 
$$
\|T(t)x\| \ge c e^{-\a t} \|x\| \qquad (x \in X, \,t>0).
$$
For each $\o > \a$, there exists a $C_0$-group extension $S$ on a Banach space $Y \supseteq X$ and $M\ge1$ such that $\|S(t)\| = \|T(t)\|$ and $\|S(-t)\| \le Me^{\o t}$ for all $t>0$.
\endprop

\proof
By replacing $T(t)$ bt $e^{\a t}T(t)$ and $\o$ by $\o-\a$, we may assume that $\alpha=0$.  Let $\o>0$ and $\tau = \o^{-1} \log (c^{-1})$.  By \thmref{3.1.1},  there exist $Y \supseteq X$ and an isometric homomorphism $\varphi: \{T(\tau)\}' \to \mathcal{B}(Y)$ such that $\varphi(T(t))$ is an extension of $T(t)$ for all $t$, $\varphi(T(\tau))$ is invertible and $\|\left(\varphi(T(\tau))\right)^{-1}\| \le c^{-1}$.  By applying \propref{2.2.0} and changing $Y$ if necessary, we may arrange that there is a $C_0$-semigroup $S$ on $Y$ given by
$$
S(t) = \varphi(T(t)) \qquad (t \ge0).
$$
Since $S(\tau)$ is invertible, $S$ is a $C_0$-group on $Y$ extending $T$.  Moreover
$$
\|S(-\tau n)\|  \le c^{-n} = e^{\o \tau n} \qquad (n\ge1).
$$
 This implies that
$$
\|S(-t)\| \le M e^{\o t} \qquad (t\ge0),
$$
where $M\ge1$ is chosen suitably. 
\endproof

Xu and Shang \cite[Theorem 2.4]{XS}  have stated the following.

\thm{XS}
Let $T$ be a $C_0$-semigroup on a Banach space $X$, with generator $A$.  The following properties are equivalent:
\begin{enumerate}[\rm(i)]
\item $T$ satisfies lower bounds,
\item There exists an equivalent norm $\|\cdot\|_*$ on $X$ such that, for some $\o \in \R$, $-(A+\o)$ is dissipative on $(X, \|\cdot\|_*)$.
\end{enumerate}
\endthm

The proof of (i)$\Rightarrow$(ii) in \cite{XS} appears to be incomplete, as the equivalent norm is taken to be 
\eqn{XS}
\|x\|_* = \inf_{t\ge0} \|e^{\a t} T(t)x\|  \qquad (x \in X),
\endeqn
where
\eqn{left}
\|T(t)x\| \ge c e^{- \a t} \|x\|  \qquad (x \in X).
\endeqn
The authors then take $\omega=\alpha$ in (ii).  However the definition of $\|\cdot\|_*$ in \eqnref{5.3} may not satisfy the triangle inequality.   So we give here a proof that (i) implies (ii), based on \propref{3.2.13}.   This proof shows that when \eqnref{left} holds for some $c$ and $\a$, then $\o$ in (ii) can be any number with $\o>\a$.  We shall see in \exref{3.2.12} that it may not be possible to take $\o=\a$.

\begin{proof} [Proof of \thmref{XS}, \rm(i)$\Rightarrow$\rm(ii)]
Let $c$ and $\a$ be as in \eqref{left}, and $\o>\a$.  Let $S$ be a $C_0$-group extension of $T$ on  $Y \supseteq X$ as in \propref{3.2.13}.  By a standard renorming in semigroup theory \cite[Lemma II.3.10]{EN00}, there is an equivalent norm on $Y$ given by
$$
\|y\|_* = \sup \{e^{-\o t}\|S(-t)y\| :t\ge0\} \qquad (y \in Y),
$$
and $e^{-\o t}S(-t)$ is a $C_0$-semigroup of contractions on $(Y, \|\cdot\|_*)$, so its generator is dissipative in that space.   The generator is $-(B+\o)$, where $B$ is the generator of the $C_0$-semigroup $S$.  By \propref{2.0.0}, $B$ is an extension of $A$.  Hence $-(A+\o)$ is dissipative on $(X,\|\cdot\|_*)$.
\end{proof}

When $A$ generates a $C_0$-semigroup $T$, the following result gives extension properties of $T$ which characterise when there is an equivalent norm for which $-A$ is dissipative.

\thm{new}  Let $T$ be a $C_0$-semigroup on a Banach space $X$, with generator $A$.  The following are equivalent:
\begin{enumerate}[\rm(i)]
\item  There is an equivalent norm $\|\cdot\|_*$ on $X$ such that $-A$ is dissipative with respect to $(X,\|\cdot\|_*)$.
\item  There is a Banach space $Y$ in which $X$ is isomorphically embedded, and a $C_0$-group $S$ on $Y$ such that $\|S(-t)\|_{\mathcal{B}(Y)} \le 1$ and $S(t)x = T(t)x$ for all $x \in X, \, t>0$.
\item  There is a $C_0$-group extension $S$ of $T$ on a Banach space $Z \supseteq X$ and a constant $\kappa$ such that $\|S(-t)\|_{\mathcal{B}(Z)} \le \kappa$ for all $t\ge0$.
\item  There exists $\kappa$ such that
$$
\|x\| \le \kappa ( \|x_1\| + \dots + \|x_n\| )
$$
whenever  $x_1,\dots, x_n \in X$, $0 \le t_1 \le t_2 \le \dots \le t_n \le t$ and
$$
T(t) x = T(t_1) x_1 + \dots + T(t_n) x_n.
$$
\end{enumerate}
\endthm

\proof (i)$\Rightarrow$(ii):  Apply \propref{3.5.7b} in the space $(X,\|\cdot\|_*)$.

\noindent (ii)$\Rightarrow$(i):  Let $\|\cdot\|_*$ be the norm of $Y$ restricted to $X$.   The generator of $\{S(-t):t\ge0\}$ is an extension of $-A$, and it is dissipative with respect to $\|\cdot\|_Y$.

\noindent (ii)$\Rightarrow$(iii):  Take $\kappa_1,\kappa_2>0$ such that $\kappa_1 \|x\|_Y \le \|x\|_X \le \kappa_2 \|x\|_Y$ for all $x \in X$.  For $y \in Y$, let
$$
\|y\|_Z = \inf \left\{ \|x\|_X + \kappa_2 \|y-x\|_Y : x \in X\right\}.
$$
It is readily verified that $\|\cdot\|_Z$ is a norm on $Y$, equivalent to $\|\cdot\|_Y$, and $\|x\|_Z = \|x\|_X$ for all $x \in X$.  Then we can take $Z = (Y,\|\cdot\|_Z)$.

\noindent (iii)$\Rightarrow$(ii):  Let $Y= (Z, \|\cdot\|_Y)$ where
$$
\|y\|_Y = \sup \{\|S(-t)y\|_Z : t \ge 0\}.
$$
This norm is equivalent to $\|\cdot\|_Z$ on $Z$, and hence equivalent to $\|\cdot\|_X$ on $X$, and $S(-t)$ is contractive with respect to $\|\cdot\|_Y$.

\noindent (iii)$\Leftrightarrow$(iv):    This follows from \thmref{3.2.4}, taking $\gamma(t) = \kappa \ge 1$ for all $t\ge0$.
\endproof

The following example shows that  \eqnref{4.0} for $c(t)=c$ does not imply the equivalent conditions of \thmref{new}.  Consequently, \eqnref{4.0} for $c(t)=c$ does not imply that $T$ is expansive in an equivalent norm on $X$, and one cannot take $\o=\a$ in \propref{3.2.13}.

\ex{3.2.12}
Let $c \in (0,1)$ and $M>1$ be given.  There exists a $C_0$-semigroup $T$ on a Hilbert space $X$ with the following properties:
\begin{enumerate}[\rm(a)]
\item $c\|x\| \le \|T(t)x\| \le \sqrt2 M^t\|x\|$ for all $t\ge0, x\in X$,
\item If $L$ is any left-inverse semigroup for $T$ on a Banach space $Y \supseteq X$ then $\{\|L(t)\| : t\ge0\}$ is unbounded.
\end{enumerate}

We begin by presenting examples satisfying a weaker version of (b).

Given $M > 1$ and $c \in (0,1)$, take $\ep \in (0,c)$.  Take $k_\ep>0$ and $\ep'>0$ such that
$$
c^k < \ep, \quad  \frac{M^k (\ep-c^k)}{2} > \frac{Mc}{1-c}, \quad \frac{\ep-c^k}{2} < \frac{M}{1-c} \ep' < \ep - c^k.
$$
Then let
\begin{align*}
\Omega &= \cb{(x,y) \in \R_+^2 : \text{$x=0$ or $y \in k_\ep\N$}}, \\
w(x,y) &= \begin{cases} c^{\lceil y/k_\ep \rceil - 1}   & \text{if $0 \le x \le t < x+y$}, \\ c^{y/k_\ep} M^x \ep' &\text{if $x>0$}. \end{cases}
\end{align*}
Let $\mu$ be one-dimensional Hausdorff measure on $\Omega$, and let $X_\ep$ be the Hilbert space $L^2(\Omega,w,\mu)$ of all $f$ for which
$$
\|f\|^2 := \int_\Omega |f|^2 w^2 \, d\mu < \infty.
$$
Let
$$
(T_\ep(t)f)(x,y) = \begin{cases} f(x-t,y) & \text{if $0\le t<x$}, \\ f(0, y -t+x) & \text{if $0\le x \le t < x+y$}, \\ 0 &\text{otherwise}. \end{cases}
$$
Then $T_\ep$ is a $C_0$-semigroup with $\|T_\ep(t)\| \le \sqrt2 M^t$ and $\|T_\ep(t)f\| \ge c\|f\|$ for all $f \in X$ and $t\ge0$.  Take $f \in X$ which is supported on $\{0\} \times (0,k_\ep)$, with $\|f\|=1$.  Then there exist $g, f_i \in X$ with 
\begin{gather*}
T_\ep(k_\ep^2)f = g + \sum_{i=1}^{k_\ep} T_\ep(k_\ep(k_\ep-i))f_i,  \qquad \|g\| + \sum_{i=1}^{k_\ep} \|f_i\| \le \ep.
\end{gather*}
We refer the reader to \cite[Section 3.2]{Gey} for the detailed justification of these properties.  It follows that for any left-inverse $C_0$-semigroup $L$ on a Banach space $Y \supseteq X_\ep$, 
$$
\sup_{0\le t \le k_\ep^2}  \|L(t)\| \ge \ep^{-1}.
$$ 

To complete the construction, take the $\ell_2$-direct sum $X$ of the spaces $X_{n^{-1}}  \; (n > 1/c)$ as constructed above with $\ep = n^{-1}$, and $k_\ep = \lceil \log(2n)/\log(1/c) \rceil$ and let $T$ be the direct sum of the $C_0$-semigroups $T_{n^{-1}}$.  Any left-inverse semigroup for $T$ on $Y \supseteq X$ is then a left-inverse semigroup for $T_{n^{-1}}$ on $X_{n^{-1}}$ (regarded as a subspace of $X$ in the natural way), so there exist $t_n$ such that 
$$
0 \le t_n \le \left(1+\log(2n)/\log(1/c)\right)^2, \qquad \|L(t_n)\| \ge n \qquad (n \ge 1/c).
$$
This implies that
$$
\limsup_{t\to\infty} e^{-\sqrt{t}} \|L(t)\| > 0.
$$
In particular $\|L(t)\|$ is not polynomially bounded in $t$.
\endex

\sect6{Operators on Hilbert spaces}

In this section we give results for operators on Hilbert spaces, where the extensions should also be operators on  Hilbert spaces.   We start by reformulating \propref{6.0.0} in a form which is specific to Hilbert spaces.

When $X$ is a Hilbert space, we may take $Y = X \ominus \ran A$ and $p=2$ in the construction of \propref{6.0.0}.  Then $Z$ is also a Hilbert space, and the extension $B$ constructed there satisfies $\|B^{-1}\| \le c^{-1}$ and $\|B\|=\|A\|$ if $A$ is bounded (note that $\|A\|\ge c$).   If in addition $A$ is densely defined we can construct the same (up to unitary equivalence) invertible extension of $A$ via the polar decomposition of $A$, $A = PU$ where $P = (A^*A)^{1/2}$ is self-adjoint and $U$ is a partial isometry.   When $A$ is bounded below as in \eqnref{2.1}, $P$ is invertible and $U$ is a unitary operator of $X$ onto the closed subspace $\ran A$.

\prop{3.3.5}
Let $A$ be a closed, densely defined operator on a Hilbert space $X$, and assume that {\rm\eqnref{2.1}} holds
for some $c>0$. Let $A = UP$ be the polar decomposition of $A$.   Let $\tilde Z$ be the Hilbert space $X \times X$,  and $\tilde B$ be the operator on $\tilde Z$ given by
$$
\tilde B = \begin{pmatrix} A& c(I-UU^*)\\ 0 & cU^* \end{pmatrix},
$$
with domain $D(\tilde B) = D(A) \times X$.
 Then $\tilde B$ has the following properties:
\begin{enumerate}[\rm(a)]
\item $\tilde B$ is an outer extension of $A$ (when $X$ is identified with $X \times \{0\}$).
\item $\tilde B$ is invertible and $\|\tilde B^{-1}\| \le c^{-1}$.
\item If $A$ is not invertible, then $\s(\tilde B) \sub \s(A)$.
\item\label{6.1d} If $A$ generates a $C_0$-semigroup on $X$, then $\tilde B$ generates a $C_0$-semigroup on $Z$.
\item If $A$ is bounded, then $\tilde B$ is bounded and $\|\tilde B\|=\|A\|$.
\end{enumerate}
\endprop

\proof
Replacing $A$ by $c^{-1}A$ we will assume that $c=1$, which simplifies the presentation.  

It is clear that $(x,0) \in D(\tilde B)$ if and only if $x \in D(A)$ and then $\tilde B(x,0) = (Ax,0)$.   Moreover $\tilde B$ is invertible with
$$
\tilde B^{-1} = \begin{pmatrix} P^{-1}U^* & 0 \\ I - UU^* & U \end{pmatrix},
$$
and $\|\tilde B^{-1}\| \le 1$ because
\begin{align*}
\|\tilde B(x_1,x_2)\|^2 &= \|Ax_1 + (I-UU^*)x_2\|^2 + \|U^*x_2\|^2 \\
&= \|Ax_1\|^2 + \|(I-UU^*)x_2\|^2 + \|UU^*x_2\|^2\\
&= \|Ax_1\|^2 + \|x_2\|^2 \ge \|x_1\|^2 + \|x_2\|^2.
\end{align*}

Assume that $A$ is not invertible.  The assumption \eqnref{2.1} on $A$ implies that the approximate point spectrum of $A$ is contained in $\{\l\in\C: |\l|\ge1\}$.  Since the boundary of $\sigma(A)$ consists of approximate eigenvalues and $0 \in \sigma(A)$ by assumption, it follows that the closed unit disc is contained in $\sigma(A)$. Thus if $\l \in \rho(A)$ then $|\l|>1$,  so $\l \in \rho(U^*)$, $R(\l,U^*) = \sum_{j=0}^\infty \l^{-(j+1)}(U^*)^j$ and the operator
\eqn{6.1}
R(\l,\tilde B) = \begin{pmatrix} R(\l,A) & R(\l,A)(I-UU^*)R(\l,U^*) \\ 0 & R(\l,U^*) \end{pmatrix}
\endeqn
is a two-sided inverse of $\l-\tilde B$.  Thus $\l \in \rho(\tilde B)$.  So $\s(\tilde B) \sub \s(A)$.  

The proof of (\ref{6.1d}) is essentially the same as in \propref{6.0.0}, as the operator $\tilde B$ is a bounded perturbation of 
$$
\begin{pmatrix} A & 0 \\ 0 & 0 \end{pmatrix}.
$$

If $A$ is bounded, then
$$
\|\tilde B(x_1,x_2)\|^2 = \|Ax_1\|^2 + \|x_2\|^2 \le \|A\|^2 \,\|(x_1,x_2)\|^2,
$$
since \eqnref{2.1} for $c=1$ implies that $\|A\|\ge1$.
\endproof

The extension of $A$ to $\tilde B$ in \propref{3.3.5} is not necessarily minimal.  In order to achieve the minimal extension $B$, we should replace $\tilde Z$ by the closure of $\bigcup_{k\ge0} \tilde B^{-k}(X \times \{0\})$.  This space is identified in the following way.

\prop{3.3.5b}
Let $X$, $\tilde Z$, $A$, $c$ and $\tilde B$ be as in \propref{3.3.5}, and identify $X$ with the subspace $X \times \{0\}$ of $\tilde Z$.  Let $Z$ be the closure of $\bigcup_{k\ge0} \tilde B^{-k}(X)$ in $\tilde Z$.  

For $k\ge0$, let 
$$
\tilde Z_k = \ran \left(U^k(I-UU^*)\right), \qquad  Z_k = \{0\} \times \tilde Z_k.
$$
The subspaces $\tilde Z_k$ of $X$ are closed and pairwise orthogonal.  Moreover
\eqn{6.2}
Z = X \oplus \bigoplus_{k=0}^\infty Z_k,
\endeqn
as an orthogonal direct sum.

Let $B$ be the restriction of $\tilde B$ to $D(B) := D(A) \oplus \bigoplus_{k=0}^\infty Z_k$, considered as an operator on $Z$.  The following properties hold:
\begin{enumerate}[\rm(1)]
\item\label{6.2(1)}  $B$ satisfies the properties {\rm(a)--(f)} of \propref{6.0.0}, when $X$ is identified with $X \times \{0\}$, and $\|B^{-1}\|\le c^{-1}$. 
\item\label{6.2(2)} $c^{-1}B$ is a unitary operator of $Z_k$ onto $Z_{k-1}$ for each $k\ge0$, where $Z_{-1} = (X \ominus \ran A) \sub X$.
\item\label{6.2(3)} Let $W \supseteq X$ be a Hilbert space with closed subspaces $W_k$ such that 
$$
W = X \oplus \bigoplus_{k=0}^\infty W_k,
$$
as an orthogonal direct sum, and let $C$ be an invertible operator on $W$ such that $\|C^{-1}\| \le c^{-1}$, $C$ extends $A$, $D(C) = D(A) \oplus (W \ominus X)$, and $C$ maps $W_k$ onto $W_{k-1}$ for $k\ge0$, where $W_{-1} = X \ominus \ran A$.  Then there is a contraction $\pi : Z \to W$ which acts as the identity on $X$, maps $D(B)$ into $D(C)$, and satisfies $C \circ \pi = \pi \circ B$ on $D(B)$.  If $c^{-1}C : W_k \to W_{k-1}$ is unitary for each $k\ge0$, then $\pi : Z \to W$ is unitary.
\end{enumerate}
\endprop

\proof  We again assume that $c=1$.

It is elementary that the spaces $\tilde Z_k$ are closed and orthogonal in $X$ and so the spaces $Z_k$ are closed and orthogonal in $\tilde Z$.  It is easily seen that $\tilde B$ maps $Z_k$ isometrically onto $Z_{k-1}$ for $k \ge 0$, and $\tilde B^{-k}(X \times \{0\}) = X \times \tilde Z_{k-1}$, for all $k\ge1$.   This establishes that \eqnref{6.2} and property (\ref{6.2(2)}) hold.  Moreover $Z$ is invariant under $\tilde B$ and $\tilde B^{-1}$, and under $R(\l,\tilde B)$ for $|\l|>1$.  So properties (a)-(e) established in \propref{3.3.5} for $\tilde B$ transfer to $B$.  To see (\ref{3.1f}), one may use (\ref{6.2(3)}), proved in the next paragraph, to deduce that the extension $(Z,B)$ is unitarily equivalent to the extension $(Z,B)$ obtained from \propref{6.0.0} with $p=2$ and $Y = Z_{-1}$, so (\ref{3.1f}) transfers from the corresponding statement in \propref{6.0.0}.  Alternatively one may carry out the calculations directly in the space $Z$ (see \cite{Gey}). 

Let $W$, $W_k$ and $C$ be as in (\ref{6.2(3)}).  For each $k\ge0$ the map $C^{-(k+1)} B^{k+1}$ is a contraction of $Z_k$ onto $W_k$.  Together with the identity map on $X$ these maps define a contraction $\pi$ of $Z$ into $W$.  The property that $C \circ \pi = \pi \circ B$ is easily seen to hold on $X$ and each $Z_k$. If $C : W_k \to W_{k-1}$ is isometric and surjective for each $k$, then so is $\pi$.  So (\ref{6.2(3)}) is established.
\endproof

As in the case of Banach spaces, it would be most useful to have a Hilbert space construction in which a given operator $A$ is extended to an invertible operator with preservation of its  spectrum and of norms of associated bounded operators, and also including the existence of a suitable homomorphism $\varphi$ from $\{A\}'$.  This would be similar to \thmref{3.1.1} but with $X$ and $Y$ both Hilbert spaces.    By adapting Arens's method, Badea and M\"uller \cite[Corollary 4.8]{BM} obtained the following result for bounded operators.  It has most of the desired features but the estimates for the norms of the powers of the inverse are not as sharp as one might expect. 

\thm{3.1.4}
Let $U$ be a bounded operator on a Hilbert space $X$ and assume that
$$
\|Ux\| \ge c \|x\|  \qquad (x \in X)
$$
for some $c>0$.  Then there exists a Hilbert space $Y \supseteq X$ and a unital homomorphism $\varphi: \{U\}' \to \mathcal{B}(Y)$ such that 
\begin{enumerate}[\rm(a)]
\item $\varphi(U)$ is invertible and $\|\varphi(U)^{-k}\| \le c^{-k} 2^{(k+1)/2}$ for all $k\ge 1$,
\item  For each $V \in \{U\}'$, $\varphi(V)$ is an extension of $V$ and $\|\varphi(V)\| \le \sqrt2 \|V\|$.
\end{enumerate}  
\endthm

We note also that Read \cite{Rea87} gave a version of \thmref{3.6.1} for bounded operators on Hilbert space, in which the norm of the extensions was almost preserved. 

\quest{3.1.6}  For which bounded operators $U$ on a Hilbert space satisfying \eqnref{2.1} do there exist an invertible extension $V$ of $U$ on a Hilbert space $Y \supseteq X$  and an isometric unital homomorphism $\varphi :  \{U\}' \to \mathcal{B}(Y)$ with $\varphi(U) = V$?
\endquest

\sect7{Semigroups on Hilbert spaces}

Recall that a densely defined operator $A$ on a Hilbert space $X$ generates a $C_0$-semigroup of contractions if and only if $A$ is a maximal dissipative operator, i.e., $A$ is dissipative and it has no proper extension which is a dissipative operator on $X$ \cite[Theorem 10.4.2]{Dav}.  Moreover any dissipative operator on $X$ has a maximal dissipative extension on $X$ (see \cite[Theorem 3.1.2]{Gor} for a simple proof).

The following result gives a different type of extension, by showing that any dissipative operator on a Hilbert space has an \textit{outer} extension which generates a contraction semigroup on a larger Hilbert space.

\thm{3.5.3}  Let $A$ be a closed, densely defined, dissipative operator on a Hilbert space $X$, and assume that there exists $\mu \in \C$ such that $\Re\mu<0$ and $\mu \notin \sp(A)$.  Then there exist a Hilbert space $Y \supseteq X$ and an outer extension $B$ of $A$ such that $B$ generates a $C_0$-semigroup of contractions on $Y$. 
\endthm

\proof
Substituting $A$ by $aA + ib$ for some $a>0$ and $b \in \R$ if necessary, we may assume that $-1 \notin \sp(A)$.  Consider the Cayley transform of $-A$:
$$
C := (I-A)(I+A)^{-1}, \qquad D(C) = \ran(I+A).
$$
Then $C$ is closed, $I+C = 2(I+A)^{-1}$ is injective with dense range $D(A)$, and $A = (I-C)(I+C)^{-1}$ on $D(A)$.  Moreover $\|Cx\| \ge \|x\|$ for all $x \in X$, by putting $y = (I+A)^{-1}x$ in
$$
\|(I-A)y\|^2 - \|(I+A)y\|^2 = - 4 \Re \langle Ay,y \rangle \ge 0.
$$
We apply \propref{3.3.5b} to the operator $C$, putting $c=1$ and $\l=-1$ in the statement of \propref{6.0.0}(f). Then $C$  has an extension to a bounded invertible operator $G$ with $\|G^{-1}\| \le 1$ on some Hilbert space $Y \supseteq X$, $\ker(I+G) = \{0\}$ and $\ran(I+G)$ is dense in $Y$, and $\ran(I+G) \cap X = \ran(I+C)$.  Thus the operator
$$
B := (I-G)(I+G)^{-1}, \qquad D(B) = \ran(I+G),
$$
is well-defined with dense domain.  Moreover
$$
I+B = 2(I+G)^{-1},  \qquad I-B = 2G(I+G)^{-1},
$$
so $1 \in \rho(B)$, $(I-B)^{-1} = \frac12 (G^{-1} + I)$, and the Cayley transform of $B$ is
$$
(I+B)(I-B)^{-1}  = G^{-1}.
$$
Since the Cayley transform of $B$ is contractive, $B$ is dissipative.  Thus $B$ generates a $C_0$-semigroup of contractions on $Y$, by the Lumer-Phillips theorem.

To show that $B$ extends $A$, we have
$$
D(B) \cap X = \ran(I + G) \cap X = \ran(I+C) = D(A).
$$
For $x \in D(A)$, $x = (I+C)x'$ for some $x' \in X$ and
\begin{multline*}
(I+G)(x+Bx) = (I+G)(2x') = 2x = (I+C)(2x') \\
= (I+C)(x+Ax) = (I+G)(x+Ax).
\end{multline*}
Since $I+G$ is injective, $Bx=Ax$.
\endproof

Now we turn to lower bounds for semigroups on Hilbert space.  The following is a Hilbert space version of \propref{3.5.7b}, in which the left inverses act on the same space as the original operator. 

\prop{3.5.7}  Let $A$ be the generator of a $C_0$-semigroup $T$ on a Hilbert space $X$.  The following are equivalent:
\begin{enumerate}[\rm(i)]
\item\label{7.2i} $T$ is expansive, i.e., $\|T(t)x\| \ge \|x\|$ for all $x \in X, t\ge0$.
\item\label{7.2ii} $-A$ is dissipative.
\item\label{7.2iii}  There is a left-inverse semigroup for $T$ consisting of contractions on $X$.
\end{enumerate}
\endprop

\proof   The equivalence of (\ref{7.2i}) and (\ref{7.2ii}) was established in \propref{3.5.7b}, and it is elementary that (\ref{7.2iii}) implies (\ref{7.2i}).

Assume that $-A$ is dissipative.  Then $-A$ has a maximal dissipative extension $B$ which generates a $C_0$-semigroup $L$ of contractions on $X$.  By \propref{2.0.0}(b), $L$ is a left-inverse semigroup for $T$.    
\endproof

We do not know whether the properties in \propref{3.5.7} are also equivalent to the existence of a $C_0$-group extension $S$ of $T$ on a Hilbert space $Y \supseteq X$ such that $\|S(-t)\| \le 1$ for all $t\ge0$.  

The following includes results of Louis and Wexler \cite[Corollary, p.260]{LW} and Zwart \cite{Zwa}.  

\thm{3.2.14}
Let $T$ be a $C_0$-semigroup on a Hilbert space $X$, with generator $A$.  The following are equivalent:
\begin{enumerate}[\rm(i)]
\item\label{7.3i}  $T$ satisfies lower bounds.
\item\label{7.3ii}  There is a left-inverse semigroup for $T$ on $X$.
\item\label{7.3iii}  There is a $C_0$-group extension $S$ of $T$ on a Hilbert space $Y \supseteq X$.
\item\label{7.3iv}  There exists $Q \in \mathcal{B}(X)$ and an equivalent inner product on $X$ such that the $C_0$-semigroup generated by $A+Q$ is isometric in the equivalent norm.
\end{enumerate}
If $c>0$ and $\a\in\R$ are such that $\|T(t)x\| \ge ce^{-\a t} \|x\|$ for all $x \in X$ and $t>0$, and $\o>\a$, then $Y$ and $S$ in {\rm(\ref{7.3iii})} may be chosen so that $\|S(t)\| \le \sqrt2 \|T(t)\|$ and $\|S(-t)\| \le M e^{\o t}$ for all $t>0$, for some constant $M$.
\endthm

\proof (\ref{7.3i})$\Rightarrow$(\ref{7.3ii}): This was proved in \cite[Corollary, p.260]{LW}, and an alternative proof was given in \cite[Theorem 1]{Zwa}.

\noindent (\ref{7.3i})$\Rightarrow$(\ref{7.3iii}):  This can be proved in a very similar way to \propref{3.2.13}, using \thmref{3.1.4} instead of \thmref{3.1.1}.

\noindent (\ref{7.3i})$\Rightarrow$(\ref{7.3iv}):  This was proved in \cite[Theorem 3]{Zwa}. 

\noindent (\ref{7.3iv})$\Rightarrow$(\ref{7.3iii}):  Let $\tilde X$ be the space $X$ with the equivalent inner product and $\tilde T$ be the isometric semigroup generated by $A+Q$ on $\tilde X$.  By a  result originally due to Cooper \cite{Coo} (see also \cite{Dou}),  $\tilde T$ has an extension to a unitary group $\tilde S$ on a Hilbert space $\tilde Y \supseteq \tilde X$.  Let $\|\cdot\|_{\tilde Y}$ be the norm on $\tilde Y$ and $P$ be the orthogonal projection of $Y$ onto $\tilde X$.  Define 
$$
\|y\|_Y^2 = \|Py\|_X^2 + \|(I-P)y\|_{\tilde Y}^2  \qquad (y \in \tilde Y).
$$
This is an equivalent norm on $\tilde Y$ and it coincides with $\|\cdot\|_X$ on $X$.  The space $Y := (\tilde Y, \|\cdot\|_Y)$ is a Hilbert space in the norm $\|\cdot\|_Y$.  Let $\tilde B$ be the generator of $\tilde S$ and let $S$ be the $C_0$-group on $Y$ generated by the bounded perturbation $B := \tilde B - QP$.  Then $B$ extends $A$ and $S$ extends $T$.

As in \sectref5 either of (\ref{7.3ii}) or (\ref{7.3iii}) implies (\ref{7.3i}).
\endproof

The proofs in \cite{LW} and \cite{Zwa} use ideas from control theory.  In contrast to \propref{3.2.13}, these arguments give little information about the norm of $\|S(t)\|$ for either $t>0$ or $t<0$.

We conclude with some remarks concerning the relations between Zwart's result \cite[Theorem 3]{Zwa}  (the implication (\ref{7.3i})$\implies$(\ref{7.3iv}) in \thmref{3.2.14}) and a similar result of Haase \cite[Theorem 3.1]{Ha04}, \cite[Theorem 7.2.8]{Haa}.   Haase's result says that if $B$ is the generator of a $C_0$-group on a Hilbert space $Y$ then $B$ is a bounded perturbation of the generator of a group of unitaries with respect to an equivalent inner product on $Y$.  In an earlier version of this paper, we raised the question whether the two results are logically related, and we are grateful to Abraham Ng for providing the key to the following arguments.

First, we assume Zwart's result, and we will show that Haase's result follows.   Let $B$ generate a $C_0$-group on a Hilbert space $X$.  Then the semigroup generated by $B$ satisfies lower bounds, so Zwart's result implies that there exist $Q \in \B(X)$ and an equivalent inner product on $X$ such that $B-Q$ generates a semigroup of isometries in the new norm.   In addition, $B-Q$ generates a $C_0$-group, so it generates a $C_0$-group of unitaries in the new norm.   This gives Haase's result.

Now, we assume Haase's result, and we will use it to show the implication (\ref{7.3iii})$\implies$(\ref{7.3iv}) in \thmref{3.2.14}.  Let $A$ be the generator of a $C_0$-semigroup $T$ on a Hilbert space $X$, with a $C_0$-group extension $S$ on a Hilbert space $Y \supseteq X$, with generator $B$.   By Haase's result, there exists $Q \in \B(Y)$ such that $B-Q$ generates a $C_0$-group of unitaries for an equivalent inner product on $Y$.   Let $\tilde Y$ and $\tilde X$ be the spaces $Y$ and $X$ with the equivalent inner product, and let $P$ be the orthogonal projection of $\tilde Y$ onto $\tilde X$.  For $x \in D(A) = D(B) \cap X$,
$$
Ax = PAx = P(B-Q)x + PQx.
$$
Let $Q_1 \in \B(X)$ be the restriction of $PQ$ to $X$, and let $A_1$ be the restriction of $P(B-Q)$ to the domain $D(A)$.  Then $A_1 =A-Q_1$ is the generator of a $C_0$-semigroup on $X$.  Since $B-Q$ is skew-symmetric on $\tilde Y$, $A_1$ is skew-symmetric on $\tilde X$.   Hence the semigroup generated by $A_1$ is a semigroup of isometries on $\tilde X$.  This establishes (\ref{7.3iv}).   

Finally, suppose that $\o > \a$ as in \thmref{3.2.14} and in addition $\o$ is greater than the exponential growth bound of $T$.  Then the equivalent inner product and the bounded operator $Q$ in (\ref{7.3iv}) can be chosen so that $Q$ is self-adjoint on $\tilde X$ and $\|Q\|_{\B(\tilde X)} \le \o$.  This follows from the corresponding statements in Haase's result combined with the proof in the paragraph above.

\end{document}